\documentclass[10pt,twoside]{article}
\usepackage{graphicx}
\usepackage{amsmath,amssymb}
\usepackage{Latex-document}

\newtheorem{theorem}{Theorem}[section]
\newtheorem{lemma}[theorem]{Lemma}

\newcommand{\Ext}{\mathrm{Ext}}
\newcommand{\Tr}{\mathrm{Tr}}

\newcommand{\im}{\mathrm{Im}}
\newcommand{\tr}{\mathrm{tr}}
\newcommand{\RE}{\mathrm{Re}}

\newcommand{\cF}{{\cal F}}
\newcommand{\bR}{{\mathbb R}}
\newcommand{\bE}{{\mathbb E}}
\newcommand{\bC}{{\mathbb C}}
\newcommand{\bN}{{\mathbb N}}
\newtheorem{definition}[theorem]{Definition}
\newtheorem{theorem/def}[theorem]{Theorem/Definition}

\newtheorem{corollary}[theorem]{Corollary}
\newtheorem{example}[theorem]{Example}
\newtheorem{problem}[theorem]{Problem}
\newtheorem{remark}[theorem]{Remark}
\newtheorem{remarks}[theorem]{Remarks}

\def\volumeno{I}
\title{\bf  Random Matrices, Free Probability\vskip -2mm
and the Invariant Subspace Problem\vskip -2mm Relative to a von Neumann
Algebra\vskip 5mm}
\author{U.~Haagerup\vspace*{-0.5cm}\thanks{Department of Mathematics \&
Computer Science, University of Southern Denmark, Campusvej 55, DK--5230
Odense M, Denmark. E-mail: haagerup@imada.sdu.dk}}
\date{\vspace{-8mm}}

\begin{document}

\maketitle

\thispagestyle{first} \setcounter{page}{273}

\begin{quotation}

\end{quotation}

\vskip 12mm

\section{Introduction} \label{section 1}\setzero
\vskip-5mm \hspace{5mm }

Random matrices have their roots in multivariate analysis in
statistics, and since Wigner's pioneering work \cite{Wi} in 1955,
they have been a very important tool in mathematical physics. In
functional analysis, random matrices and random structures have in
the last two decades been used to construct Banach spaces with
surprising properties. After Voiculescu in 1990--1991 used random
matrices to classification problems for von Neumann algebras, they
have played a key role in von Neumann algebra theory (cf.
\cite{V7}). In this lecture we will discuss some new applications
of random matrices to operator algebra theory, namely applications
to classification problems for $C^*$-algebras and to the invariant
subspace problem relative to a von Neumann algebra.

The rest of this lecture is divided into eight sections:

\begin{itemize}
\item[2.] Selfadjoint random matrices and Wigner's semicircle law.
\item[3.] Free probability and Voiculescu's random matrix model.
\item[4.] $\Ext(C_r^*(F_k))$ is not a group for $k\ge 2$.
\item[5.] Other applications of random matrices to $C^*$-algebras.
\item[6.] The invariant subspace problem relative to a von Neumann
algebra.
\item[7.] The Fuglede-Kadison determinant and Brown's spectral
distribution measure.
\item[8.] Spectral subspaces for operators in II$_1$-factors.
\item[9.] Voiculescu's circular operator $Y$ and the strictly upper triangular
operator $T$.
\end{itemize}

\section{Selfadjoint random matrices and Wigner's semicircle law} \label{section 2}
\setzero\vskip-5mm \hspace{5mm }

\setcounter{equation}{0}

A random matrix $X$ is an $n\times n$ matrix whose entries are real or
complex random variables on a probability space $(\Omega,\cF,P)$. We
denote by SGRM$(n,\sigma^2)$ the class of selfadjoint random matrices
\[
X_n = (X^{(n)}_{ij})^n_{i,j=1}
\]
where $X_{ij}$, $i,j=1,\dots,n$ are $n^2$ complex random variables and
\[
(X_{ii}^{(n)})_i,\quad (\sqrt{2}\, \RE\, X_{ij}^{(n)})_{i<j},\quad
(\sqrt{2}\, \im\,X_{ij}^{(n)})_{i<j}
\]
are $n^2$ independent identical distributed real Gaussian random
variables with mean value 0 and variance $\sigma^2$. In the
terminology of Mehta's book \cite{Me}, $X_n$ is a Gaussian unitary
ensemble (GUE). In the following we put $\sigma^2=\frac1n$ which
is the normalization used in Voiculescu's random matrix paper
\cite{V3}. By results of Gaudin, Mehta and Wigner from 1960--1965,
the joint distribution of the eigenvalues (in random order) of $X$
has density $g$ given by
\[
g_n(\lambda_1,\dots,\lambda_n)=
c_n\prod_{i<j}(\lambda_j-\lambda_i)^2\exp\big(-\frac{n}{2}\sum^n_{i=1}\lambda^2_i\big)
\]
where $c_n$ is a normalization constant, and the (average) density for a
single eigenvalue is given by
\[
h_n(x)=\frac{1}{\sqrt{2n}} \sum^{n-1}_{k=0}
\varphi_k\big(\sqrt{\frac{n}{2}}x\big)^2
\]
where $\varphi_0,\varphi_1,\dots$ is the sequence of Hermite
functions. Moreover,
\[
\lim_{n\to\infty} h_n(x)=\frac{1}{2\pi}\sqrt{4-x^2}\, 1_{[-2,2]}(x),\quad
x\in\bR
\]
(cf.\ \cite{Me}). This is Wigner's semicircle law for the GUE-case. In
the sense of weak convergence of probability measures, the semicircle
law can be proved under much more general assumptions on the entries
(see Wigner \cite{Wi}). Arnold proved in 1967 that the corresponding
strong law also holds, i.e. for almost all $\omega$ in the probability
space $\Omega$, the empirical eigenvalue distribution of $X_n(\omega)$
converges weakly to the semicircular distribution
$\frac{1}{2\pi}\sqrt{4-x^2}\,1_{[-2,2]}(x)dx$ as $n\to\infty$. Very
interesting research have been carried out on the level spacing of the
eigenvalues in the bulk of the spectrum (cf.\ \cite{Me}) and more recently
near the boundary of the spectrum (cf.\ \cite{TW1}, \cite{TW2}) for
selfadjoint Gaussian random matrices with real, complex or symplectic
entries (the GOE, GUE and GSE cases), but this is outside the scope of
the present lecture.

\section{Free probability and Voiculescu's random matrix model} \label{section 3} \setzero\vskip-5mm \hspace{5mm }

\setcounter{equation}{0}

Voiculescu proved in 1991 \cite{V3} an extensive generalization of
Wigner's semicircle law to families of independent random matrices. In
order to state the result, we will need some basic concepts from free
probability theory (cf.\ \cite{V1}, \cite{V2} and \cite{VDN}).\vskip
-1.5mm

\begin{definition}
\label{def3-1}
{\rm \cite{V1}}

\begin{enumerate}
\item A non-commutative probability space is a pair $(A,\varphi)$
consisting of a
unital complex algebra $A$ and a functional $\varphi\colon A\to\bC$ such
that $\varphi(1_A)=1$.
\item A $C^*$-probability space is a pair $(A,\varphi)$ consisting of a unital
$C^*$-algebra $A$ and a state $\varphi\colon A\to\bC$ on $A$.
\end{enumerate}
\end{definition}\vspace*{-1.5mm}

The connection to classical probability theory on a probability space
$(\Omega,\cF,P)$ is obtained by putting
\[
A = \bigcap^\infty_{p=1} L^p(\Omega)
\]
and
\[
\varphi(a) = \bE(a)=\int_\Omega a(\omega)dP(\omega),\quad a\in A
\]
or $A'=L^\infty(\Omega,P)$ with the same definition of $\varphi$. The
latter example is a $C^*$-probability space. To fit random
matrices (of size $n$) into this framework, one must instead consider the non-commutative
algebra
\[
A_n = \bigcap^\infty_{p=1} L^p(\Omega,M_n(\bC))
\]
with functional
\[
\varphi_n(a) = \bE(\tr_n(a))=\int_\Omega \tr_n(a(\omega))d\omega
\]
where $\tr_n=\frac1n \Tr$ is the normalized trace on $M_n(\bC)$.\vskip
-1.5mm

\begin{definition}\label{def3-2} {\rm \cite{V1}, \cite{V2}}
\begin{enumerate}
\item A family $(a_i)_{i\in 1}$ of elements in a non-commutative
probability space is a {\em free family} if for all $n\in\bN$ and all
polynomials $p_1,\dots,p_n\in\bC[X]$, one has
\[
\varphi(p_1(a_{i_1})\cdot\ldots\cdot p_n(a_{i_n}))=0
\]
whenever $i_1\ne i_2\ne\dots\ne i_n$ (neighbouring indices are different)
and \newline $\varphi(p_k(a_{i_k}))=0$ for $k=1,\dots,n$.
\item A family $(x_i)_{i\in j}$ of elements in a $C^*$-probability space
$(A,\varphi)$ is called a semicircular family if $(x_i)_{i\in I}$ is a
free family, $x_i=x_i^*$, $\varphi(x_i^{2k-1})=0$  and
\[
\varphi(x_i^{2k}) = \frac{1}{2\pi} \int^2_{-2}
t^k\sqrt{4-t^2}dt=\frac{1}{k+1} \begin{pmatrix} 2k\\k\end{pmatrix}
\]
for all $k\in\bN$ and all $i\in I$.
\end{enumerate}
\end{definition}\vspace*{-1.5mm}

We can now formulate Voiculescu's generalization of Wigner's semicircle
law: \vspace*{-1.5mm}

\begin{theorem}
\label{thm3.3} {\rm \cite{V3}}
Let $I$ be an index set and let for each $n\in\bN$, $(X_i^{(n)})_{i\in
I}$ be a family of independent SGRM$(n,\frac1n)$-distributed selfadjoint
random matrices. Then asymptotically as $n\to\infty$ $(X_i^{(n)})_{i\in
I}$ is a semicircular family, i.e. if $(x_i)_{i\in I}$ is a semicircular
family index by $I$ in a $C^*$-probability space
$(A,\varphi)$ then
\begin{equation}
\label{eq3-1}
\lim_{n\to\infty} \bE\, \tr_n(X_{i_1}^{(n)}\cdot\ldots\cdot
X_{i_p}^{(n)})=\varphi(x_{i_1}\cdot\ldots\cdot x_{i_p})
\end{equation}
for all $p\in\bN$ and all $i_1,\dots,i_p\in I$.
\end{theorem}\vskip -4mm

The corresponding strong law: For almost all $\omega\in\Omega$, one has
\begin{equation}
\label{eq3-2} \lim_{n\to\infty}
\tr_n(X_{i_1}^{(n)}(\omega)\cdot\ldots\cdot
X_{i_p}^{(n)}(\omega))=\varphi(x_{i_1}\cdot\ldots\cdot x_{i_p}),
\end{equation}
whick was proved independently by Hiai and Petz \cite{HP2} and
Thorbjørnsen \cite{T}.

\section{\boldmath$\Ext(C^*_r(F_k))$ is not a group for $k\ge 2$}
\label{section 4}  \setzero\vskip-5mm \hspace{5mm }

\setcounter{equation}{0}

Very recently Thorbjørnsen and the lecturer proved that the strong
version \eqref{eq3-2} of Voiculescu's random matrix model also holds for
the operator norm:\vspace*{-1.5mm}

\begin{theorem}
\label{thm4-1} {\rm \cite{HT3}}
Let $r\in\bN$ and let for each $n\in\bN$
$(X_1^{(n)},\dots,X_r^{(n)})$ be a set of $r$ independent
SGRM$(n,\frac1n)$-distributed selfadjoint random matrices. Let further
$(x_1,\dots,x_r)$ be a semicircular system in a $C^*$-probability space
$(A,\varphi)$, where $\varphi$ is a faithful state on $A$. Then there is
a null set $N\subseteq\Omega$ such that for all
$\omega\in\Omega\backslash N$ and all non-commutative polynomials $P$ in
$r$ variables
\[
\lim_{n\to\infty} \|P(X_1^{(n)}(\omega),\dots, X_r^{(n)}(\omega))\| =
\|P(x_1,\dots,x_r)\|.
\]
\end{theorem}\vspace*{-1.5mm}

Let $\Gamma$ be a countable (discrete) group. The reduced group
$C^*$-algebra $C^*_r(\Gamma)$ is the $C^*$-subalgebra of
$B(\ell^2(\Gamma))$ generated by the set of unitaries
$\{\lambda(\gamma)\mid \gamma\in\Gamma\}$, where $\lambda\colon\Gamma\to
B(\ell^2(\Gamma))$ is the left regular representation. By the methods of
\cite{V2} it follows that for the free group $F_k$ on $k$ generators,
$C_r^*(F_k)$ can be embedded in $C^*(x_1,\dots,x_k,1)$, where
$x_1,\dots,x_k$ is a free semicircular family in a $C^*$-probability
space $(A,\varphi)$ with $\varphi$ faithful. Hence as a corollary of
Theorem \ref{thm4-1} we have\vspace*{-1.5mm}

\begin{corollary}
\label{cor4-2} {\rm \cite{HT3}} czj Let $k \in\bN$, $k\ge 2$. Then
$C^*_r(F_k)$ can be embedded in the quotient $C^*$-algebra $\prod
M_n(\bC)/\sum M_n(\bC)$ where
\begin{eqnarray*}
\prod M_n(\bC)&=& \left\{ (x_n)^\infty_{n=1}\mid x_n\in M_n(\bC),\
\sup_n\|x_n\|<\infty \right\}\\
\sum M_n(\bC)&=& \left\{ (x_n)^\infty_{n=1}\mid x_n\in M_n(\bC),\
\lim_{n\to\infty}\|x_n\|= 0 \right\} .
\end{eqnarray*}
In particular $C_r^*(F_k)$ is a MF-algebra in the sense of Blackadar and
Kirchberg {\rm \cite{BK}}.
\end{corollary}

The invariant $\Ext(A)$ for a $C^*$-algebra $A$ was introduced by Brown,
Douglas and Fillmore in \cite{BDF}. $\Ext(A)$ is the set of all
essential extensions $B$ of $A$ by the compact operators $K$ on the
Hilbert space $\ell^2(\bN)$, and it has a natural semigroup structure.
Voiculescu proved in \cite{V0} that $\Ext(A)$ is always a unital
semigroup, and by Choi and Effros \cite{CE} $\Ext(A)$ is a group, when
$A$ is a nuclear $C^*$-algebra. Andersen \cite{An} provided in 1978 the
first example of a $C^*$-algebra $A$ for which $\Ext(A)$ is not a group.
The $C^*$-algebra in \cite{An} is generated by $C^*_r(F_2)$ and a
projection $p\in B(\ell^2(F_2))$. Since then it has been an open problem
whether $\Ext(C^*_r(F_2))$ is a group (see \cite[Sect.5]{V5} for a more
detailed discussion about this problem). It is well known that a proof
of Corollary \ref{cor4-2} would provide a negative solution to this
problem (see \cite[5.12]{V5}, \cite{V4} and \cite{Ro}). The argument
works for all $k\ge 2$. Hence we have\vspace*{-1.5mm}

\begin{corollary}
\label{cor4-3} {\rm \cite{HT3}}
For all $k\in \bN$, $k\ge 2$, $\Ext(C^*_r(F_k))$ is not a group.
\end{corollary}

\begin{remarks}
\label{rem4-2}  \ \\ \rm
 a) Corollaries \ref{cor4-2} and
\ref{cor4-3} also hold for $k=\infty$.

\noindent b) $C^*_r(F_k)$ is not quasidiagonal (cf {\rm \cite{Ro}}) but the
non-invertible extension $B$ of $C^*_r(F_k)$ obtained from Corollary
\ref{cor4-2} is quasidiagonal.

\noindent c) $C^*_r(F_k)$ is an exact $C^*$-algebra, but for any
non-invertible extension $B$ of $C^*_r(F_k)$ by the compact operators,
$B$ cannot be exact. This follows from the Lifting theorem in
{\rm \cite{EH}}. Other examples of non-exact extensions of exact
$C^*$-algebras by $K$ are given in {\rm \cite{Ki2}}.
\end{remarks}

In the rest of this section, I will briefly outline the main steps in
the proof of Theorem
\ref{thm4-1}. From \eqref{eq3-2} it follows that for all non-commutative
polynomials $P$ in $r$ variables
\begin{equation}
\label{eq5-1}
\liminf_{n\to\infty} \|P(X_1^{(n)}(\omega),\dots,X_r^{(n)}(\omega))\|\ge
\|P(x_1,\dots,x_r)\|
\end{equation}
for almost all $\omega\in\Omega$ (see \cite{T}), so we ``only'' have to
prove that
\begin{equation}
\label{eq5-2}
\limsup_{n\to\infty} \|P(X_1^{(n)}(\omega),\dots,X_r^{(n)}(\omega)\|\le
\|P(x_1,\dots,x_r)\|
\end{equation}
for almost all $\omega\in\Omega$. Even the case $r=1$ and $P(x)=x$ is a
difficult task. It corresponds to proving that if $X_n$ is
SGRM$(n,\frac1n)$-distributed, $n=1,2,\dots$ then for almost all
$\omega\in\Omega$,
$$ \limsup_{n\to\infty}\lambda_{\mbox{max}}(X_n(\omega)) \le 2 \qquad
\liminf_{n\to\infty}\lambda_{\mbox{min}}(X_n(\omega)) \ge -2, $$
where
$\lambda_{\mbox{max}}$ and $\lambda_{\mbox{min}}$ are the smallest and
largest eigenvalue of $X_n(\omega)$. This problem was settled by Bai and
Yin \cite{BY} in 1988 using Geman's combinatorial method \cite{Ge}. (See
also \cite[Thm.\ 2.12]{Ba} and \cite[Thm.\ 3.1]{HT0}).\vskip -2mm

\begin{lemma}[The linearization trick] {\rm \cite{HT3}}
\label{lemma5-1}
In order to prove \eqref{eq5-2} it is sufficient to show that for all
$m\in\bN$ and all selfadjoint $m\times m$-matrices $a_0,\dots,a_r$ and
all $\varepsilon>0$,
\begin{equation}
\label{eq5-3}
\sigma\big(a_0 \otimes 1 + \sum^r_{i=1} a_i\otimes X_i^{(n)}(\omega)\big)
\subseteq \sigma\big( a_0\otimes 1 +\sum^r_{i=1} a_i\otimes
x_i\big)+]-\varepsilon,\varepsilon[
\end{equation}
holds eventually as $n\to\infty$ for almost all $\omega\in\Omega$. Here
$\sigma(T)$ denotes the spectrum of a matrix or an operator $T$.
\end{lemma}\vspace*{-3mm}

\begin{lemma}\label{lemma5-2} {\rm \cite{HT3}}
Let $a_0,\dots,a_r$ be as above, and put
\begin{eqnarray*}
S_n &=& a_0\otimes 1 + \sum^r_{i=1} a_i\otimes X_i^{(n)}\\
s &=& a_0\otimes 1 + \sum^r_{i=1} a_i\otimes x_i.
\end{eqnarray*}
Moreover, let $G_n,G$ be the matrix valued Stieltjes transforms of $S_n$
and $S$, i.e. for $\lambda\in M_n(\bC)$, and $\mathrm{Im}\,\lambda =
\frac{1}{2i}(\lambda-\lambda^*)$ positive definite
\begin{eqnarray*}
G_n(\lambda) &=& \bE ((\mbox{id}_m\otimes\tr_n)((\lambda\otimes
1-S_n)^{-1}))\\
G(\lambda) &=& (\mbox{id}_m\otimes\varphi)((\lambda\otimes 1-s)^{-1}).
\end{eqnarray*}
Then $G_n(\lambda)$ and $G(\lambda)$ are invertible and
\begin{eqnarray}
\label{eq5-4}
a_0 + \sum^r_{i=1} a_iG(\lambda)a_i + G(\lambda)^{-1} &=&\lambda\\
\label{eq5-5}
\big\| a_0 + \sum a_iG_n(\lambda)a_i+G_n(\lambda)^{-1}-\lambda\big\| &\le &
\frac{C}{n^2}(K+\|\lambda \|)^2\|(\mbox{Im}\lambda)^{-1}\|^5
\end{eqnarray}
where $C=\frac{\pi^2m^3}{8} \big(\sum^r_{i=1} \|a_i\|^2\big)^2$ and
$K=\|a_0\|+4\sum^r_{i=1} \|a_i\|$.
\end{lemma}

The equality \eqref{eq5-4} was proved by Lehner (cf.\ \cite[Prop.4.1]{Le} using
Voiculescu's $R$-transform with amalgamation \cite{V6}. The inequality
\eqref{eq5-5} is more difficult. It relies on the concentration
phenomena used in Banach space theory, in form of \cite[Theorem
4.7]{P1}. (See \cite{Mi} for a general discussion of the concentration
phenomena.) Next we derive from \eqref{eq5-4} and \eqref{eq5-5} that
\begin{equation}
\label{eq5-6}
\| G_n(\lambda)-G(\lambda)\|\le\frac{4C}{n^2} (K+\|\lambda\|)^2
\|(\mbox{Im}\, \lambda)^{-1}\|^7
\end{equation}
when $\lambda\in M_m(\bC)$ and $\mbox{Im}\, \lambda$ is positive
definite. The estimate \eqref{eq5-6} implies that for every $f\in
C^\infty_c(\bR)$
\begin{equation}
\label{eq5-7}
\bE((\tr_m\otimes\tr_n)(f(S_n))) =
(\tr_m\otimes\varphi)(f(s))+O\big(\frac{1}{n^2}\big)
\end{equation}
for $n\to\infty$. Moreover a second application of the concentration
phenomena gives
\begin{equation}
\label{eq5-8}
\mbox{Var}((\tr_m\otimes\tr_n)(f(S_n)))
\le\frac{\pi^2}{8n^2}\bE((\tr_m\otimes\tr_n)(f'(S_n)^2))
\end{equation}
where $\mbox{Var}$ denotes the variance. Now let $g$ be a
$C^\infty(\bR)$-function with values in $[0,1]$ such that $g$ vanishes
on $\sigma(S)$ and $g$ is 1 on the complement of
$\sigma(s)+]-\varepsilon,\varepsilon[$. By applying \eqref{eq5-7} and
\eqref{eq5-8} to $f=g-1$, one gets
\begin{eqnarray}
\label{eq5-9}
\bE((\tr_m\otimes\tr_n)(g(S_n)) &=& O\big(\frac{1}{n^2}\big)\\
\label{eq5-10}
\mbox{Var}((\tr_m\otimes\tr_n)g(S_n)) &=& O\big(\frac{1}{n^4}\big).
\end{eqnarray}
By a standard application of the Borel-Cantelli lemma \eqref{eq5-9} and
\eqref{eq5-10} imply
\[
(\tr_m\otimes\tr_n)(g(S_n(\omega)))=O(n^{-4/3})
\]
almost surely. Hence the number of eigenvalues for $S_n(\omega)$ outside
$\sigma(s)+]-\varepsilon,\varepsilon[$ is $O(n^{-1/3})$\footnote{$\tr_m$
and $\tr_n$ are the normalized traces on $M_m(\bC)$ and $M_n(\bC)$.}
almost surely, but being an integer, the number has to vanish eventually
as $n\to\infty$ for almost all $\omega\in\Omega$. Hence \eqref{eq5-3}
holds.

\section{\boldmath Other applications of random matrices to $C^*$-algebras}
\label{section 6}  \setzero\vskip-5mm \hspace{5mm }

\setcounter{equation}{0}

A $C^*$-algebra $A$ is called exact if for every short exact sequence of
$C^*$-algebras
\[
0\to J\to B\to B/J\to 0
\]
the sequence
\[
0\to A\otimes_{\mbox{min}}J\to A\otimes_{\mbox{min}}B\to
A\otimes_{\mbox{min}}(B/J) \to 0
\]
is exact (cf.\ \cite{Ki1}, \cite{Wa}). The class of exact $C^*$-algebras
is very large: All nuclear $C^*$-algebras are exact and the reduced
group $C^*$-algebra $C^*_r(\Gamma)$ is exact for any discrete subgroup
$\Gamma$ of a connected locally compact group (cf.\ \cite{Ki2}). In 1991
the lecturer proved that 2-quasitraces on unital exact $C^*$-algebras
are traces (cf.\ \cite{Haa1}). Combined with results of Handelman
\cite{Han} and Blackadar and Rørdam \cite{BR}, this implies that
\begin{align}
\label{eq6-1} &\mbox{Every stably finite exact unital $C^*$-algebra has a
tracial state.}\\[7pt]
\label{eq6-2} & \mbox{Every state on the $K_0$-group, $K_0(A)$ of an exact
unital}   \\
\nonumber & \mbox{$C^*$-algebra $A$ is induced by a tracial state on $A$.}
\end{align}
Later, Thorbjørnsen and the lecturer found new proofs based on random
matrices for \eqref{eq6-1} and \eqref{eq6-2}. The key step in the proof
was to show:\vspace*{-1.5mm}

\begin{theorem} \label{thm6-1} {\rm \cite{HT1}}
Let $A$ be an exact unital $C^*$-algebra, and let $a_1,\dots,a_r\in A$
be elements in $A$ for which
\begin{eqnarray}
\label{eq6-3}
\sum^r_{i=1} a_i^*a_i &=& c\bf{1}_A\qquad\mbox{where $c>1$}\\
\label{eq6-4}
\sum^r_{i=1} a_ia_i^* &\le & \bf{1}_A
\end{eqnarray}
and let $Y_1^{(n)},\dots,Y_r^{(n)}$ be random $n\times n$-matrices whose
entries are $rn^2$ independent identically distributed complex Gaussian random
variables with density $\frac{n}{\pi}\exp(-n|z|^2)$, $z\in\bC$. Put
\begin{equation}
\label{eq6-5}
S_n = \sum^r_{i=1} a_i\otimes Y_i^{(n)}
\end{equation}
and let $\sigma(S^*_nS_n)$ be the spectrum of $S_n^*S_n$ as a function
of $\omega\in\Omega$ (the underlying probability space). Then for
almost all $\omega\in\Omega$
\begin{eqnarray}
\label{eq6-6}
\limsup_{n\to\infty} \max(\sigma(S_n^*S_n))&\le& (\sqrt{c}+1)^2\\
\label{eq6-7}
\liminf_{n\to\infty} \min(\sigma(S_n^*S_n))&\ge& (\sqrt{c}-1)^2
\end{eqnarray}
\end{theorem}\vspace*{-1.5mm}

The result is a kind of generalization of the results of Geman 1980
\cite{Ge} and Silverstein 1985 \cite{Si} on the asymptotic behaviour of
the largest and smallest eigenvalue of a random matrix of Wishart
type. The estimates \eqref{eq6-6} and \eqref{eq6-7} were proved by careful moment
estimates and lengthy combinatorial arguments. With Theorem \ref{thm4-1}
at hand, a much simpler proof of \eqref{eq6-6} and \eqref{eq6-7} can now
be obtained (cf.\ \cite{HT3}).

Theorem \ref{thm6-1} is not true in the general non-exact case
(cf.\ \cite{HT2}). It is unknown whether \eqref{eq6-1} or
\eqref{eq6-2} hold for general $C^*$-algebras. Both problems are
equivalent to Kaplansky's problem from the 1950's: Is every
AW$^*$-factor of type II$_1$ a von Neumann factor of type II$_1$?

Let me end this section by discussing another application of Theorem
\ref{thm4-1}:\\
Junge and Pisier proved in \cite{JP} that
\begin{equation}
\label{eq6-8}
B(H)\otimes_{\max}B(H)\ne B(H)\otimes_{\min}B(H).
\end{equation}
In the proof they consider a sequence of constants $C(k)$,
$k\in\bN$: For fixed $k\in\bN$ $C(k)$ is the infimum of all $C>0$ for
which there exists a sequence of $k$-tuples of unitary matrices
$(u_1^{(m)},\dots,u_k^{(m)})_{m\in\bN}$ of size $n(m)\in\bN$, such that for all
$m\ne m'$:
\[
\big\| \sum^k_{i=1} u_i^{(m)} \otimes u_i^{(m')}\|\le C.
\]
To obtain \eqref{eq6-8}, Junge and Pisier proved that
$\lim_{k\to\infty}\frac{C(k)}{k}=0$. Subsequently, Pisier
\cite{P2} proved that $C(k)\ge 2\sqrt{k-1}$ for all $k\in\bN$ and
Valette \cite{V} proved, using Ramanujan graphs, that $C(k)\le
2\sqrt{k-1}$ when $k$ is of the form $k=p+1$ for an odd prime
number $p$. It is an easy consequence of Corollary \ref{cor4-2}
that $C(k)\le 2\sqrt{k-1}$ for all $k\ge 2$ and hence
$C(k)=2\sqrt{k-1}$ for all $k\ge 2$ (see \cite{HT3}).

\section{The invariant subspace problem relative to a von Neumann algebra}
\label{section 7}  \setzero\vskip-5mm \hspace{5mm }

\setcounter{equation}{0}

The invariant subspace problem for operators on general Banach
spaces were settled by Enflo \cite{E} and Read \cite{Re} in the
1980's, but for Hilbert spaces the problem is still open:
\begin{problem} \label{prob7-1} \rm \cite[pp.~100--101]{Hal}
Let $H$ be a separable infinite dimensional Hilbert space, and let $T\in
B(H)$. Does there exist a non-trivial closed $T$-invariant subspace of $H$?
\end{problem}

More generally, one has the invariant subspace problem relative to a von
Neumann algebra:
\begin{problem}
\label{prob7-2} \rm  Let $M\subseteq B(H)$ be a von Neumann
algebra on a separable Hilbert space $H$, and let $T\in M$. Does
there exist a non-trivial closed $T$-invariant subspace $K$ for
$T$, such that $K$ is affiliated with $M$ (i.e.\ $K$ is of the
form $K=P(H)$ for a projection $P\in M$)?
\end{problem}

The problem is only interesting when $\dim(M)=+\infty$ and when $M$ is a
factor, i.e. when the
center of $M$ is just $\bC 1_M$.

The infinite dimensional factors were divided into 4 types by
Murray and von Neumann in the late 1930's (cf.\ \cite[Vol.2]{KR}).

\begin{description}
\item[Type I$_\infty$:] These are isomorphic to $B(K)$ for some infinite
dimensional Hilbert space.
\item[Type II$_1$:] $M$ has a tracial state, i.,e. there exists a
functional $\tr\colon M\to\bC$, such that $\tr(1_M)=1$, $\tr(S^*S)\ge 0$
and $\tr(ST)=\tr(TS)$ for all $S,T\in M$.
\item[Tupe II$_\infty$:] $M\simeq N\widehat{\otimes} B(K)$ where $N$ is
type II$_1$ and $\dim\, K=+\infty$.
\item[Type III:] All other infinite dimensional factors.
\end{description}
In all 4 cases, problem 2 remains open (the Type I$_\infty$ case is of
course equivalent to Problem 7.1). We will in the following address the
invariant subspace problem relative to a factor of type II$_1$.

\section{The Fuglede-Kadison determinant and Brown's spectral
distribution measure}
\label{section 8}  \setzero\vskip-5mm \hspace{5mm }

\setcounter{equation}{0}

Let $M$ be a II$_1$-factor. Then $M$ has a unique tracial state $\tr$,
and $\tr$ is normal and faithful (see eg. \cite[Vol.2, Sect.8]{KR}. The
{\em Fuglede-Kadison determinant} $\Delta\colon M\to[0,\infty)$ can be
defined (cf.\ \cite{FK}) by:
\begin{equation}
\label{eq8-1}
\Delta(T) = \lim_{\varepsilon\downarrow 0}
\exp(\tr(\log(T^*T+\varepsilon 1)^{\frac12})),\quad t\in M.
\end{equation}
If $T$ is invertible, one has
\[
\Delta(T) = \exp(\tr(\log|T|))
\]
where $|T|=(T^*T)^\frac12$. Moreover $\Delta$ has the following
properties:
\begin{eqnarray*}
\Delta(ST) &=& \Delta(S)\Delta(T),\quad S,T\in M\\
\Delta(T) &=& \Delta(T^*)=\Delta(|T|),\quad T\in M\\
\Delta(U) &=& 1,\quad\mbox{when $U\in M$ is unitary.}
\end{eqnarray*}
$\Delta$ is an upper semi-continuous function on $M$ but it is not
continuous in the norm-topology on $M$.\vspace*{-1.5mm}

\begin{theorem} \label{thm8-1} {\rm (L.G.\ Brown 1983 \cite{Br})}
Let $M$ be a II$_1$-factor and let $T\in M$. Then the function
\[
\varphi\colon\lambda\to \frac{1}{2\pi} \log\Delta (T-\lambda
1),\quad\lambda\in\bC
\]
is subharmonic and its Laplacian taken in distribution sense
\begin{equation}
\label{eq8-3}
\mu_T =
\bigg(\frac{\partial^2}{\partial\lambda^2_1}+\frac{\partial^2}{\partial\lambda_2^2}\bigg)\varphi
\end{equation}
($\lambda_1=\mbox{Re}\,\lambda$, $\lambda_2=\mbox{Im}\,\lambda$) is a
probability measure in $\bC$ concentrated on the spectrum $\sigma(T)$ of
$T$.
\end{theorem}\vspace*{-1.5mm}

\begin{definition}\label{def8-2}
The above measure $\mu_T$ is called Brown's spectral distribution
measure for $T$ or just the Brown measure for $T$.
\end{definition}

\begin{example} \label{ex8-3} \ \\ \rm
a) The Fuglede-Kadison determinant and the Brown measure also make sense
for $M=M_n(\bC)$, and $\tr=\frac1n\Tr$ the normalized trace on
$M_n(\bC)$. In this case one gets
\begin{eqnarray*}
 \Delta(T) & = & \sqrt[n]{|\det T|} \\
 \mu_T & = & \frac1n \sum^n_{i=1} \delta_{\lambda_i},
\end{eqnarray*}
where $\lambda_1,\dots,\lambda_n$ are the eigenvalues of $T$ repeated
according to root multiplicity, and $\delta_{\lambda_i}$ is the Dirac
measure at $\lambda_i$.

\noindent b) If $T$ is a normal operator (i.e. $T^*T=TT^*$) in a factor of
type II$_1$, $T$ has a spectral resolution
\[
T = \int_{\sigma(T)} \lambda dE(\lambda).
\]
In this case $\mu_T$ is equal to $\tr\circ E$.
\end{example}

Methods for computing Brown measures have been developed by Larsen and
the lecturer \cite{HL} and by Biane and Lehner \cite{BL}.

\section{\boldmath Spectral subspaces for operators in II$_1$-factors}
\label{section 9}  \setzero\vskip-5mm \hspace{5mm }

\setcounter{equation}{0}

In 1968, Apostol \cite{Ap} and Foias \cite{Fo1}, \cite{Fo2} introduced
the notion of spectral subspaces for certain well behaved operators on
Banach spaces, the decomposable operators (see \cite{LN} for a modern
treatment of this theory):\vspace*{-1.5mm}

\begin{definition} \label{def9-1} {\rm \cite[Definition 1.1.1]{LN}}
An operator $T$ on a Banach space $X$ is called decomposable if for any
open covering $\bC=V\cup W$ of the complex plane, there exist closed
$T$-invariant subspaces $Y,Z$ of $X$ such that
\begin{eqnarray}
\label{eq9-1}
X &=& Y+Z\\
\label{eq9-2}
\sigma(T|_Y) &\subseteq & V \quad\mbox{and $\sigma(T|_Z)\subseteq W$.}
\end{eqnarray}
\end{definition}

If $T\in B(X)$ is decomposable, it has a {\em spectral capacity},
i.e. there exists a map $E$ from the closed subsets of $\bC$ into the
closed $T$-invariant subspaces of $X$, such that
\begin{eqnarray}
\label{eq9-3}
E(\emptyset) & =& 0\quad\mbox{and $E(\bC)=X$}\\
\label{eq9.4}
X &=& E(\bar{V}_1)+\dots + E(\bar{V}_N)\quad\mbox{for every finite}\\ \nonumber
&& \mbox{open covering $\bC=V_1\cup V_2\cup\dots\cup V_n$}\\
\label{eq9-5}
E(\cap^\infty_{n=1} F_n)&=&\cap^\infty_{n=1} E(F_n),\quad F_n\subseteq\bC\
\mbox{closed}\\
\label{eq9-6}
\sigma(T_{|E(F)})&\subseteq & F,\quad F\subseteq\bC \  \mbox{closed}.
\end{eqnarray}
Moreover, a spectral capacity is unique (cf.\ \cite[Sect.1]{LN}).

In this section we will discuss a new method for constructing spectral
subspaces of operators which works for all operators in ``almost all''
II$_1$-factors, regardless of whether the operator is decomposable in
the above sense.\vspace*{-1.5mm}

\begin{definition}
\label{def9-2}
A II$_1$-factor $M$ on a separable Hilbert space has the {\em embedding
property} if it can be embedded in the ultrapower $R^\omega$ of the
hyperfinite II$_1$-factor $R$ for some free ultrafilter $\omega$ on the
natural numbers.
\end{definition}

All II$_1$-factors of current interest have this embedding property, and
in fact no counterexamples are known. The question whether every
II$_1$-factor on a separable Hilbert space can be embedded in $R^\omega$
was first raised by Connes in 1976 \cite{Co} (see also \cite{Ki2} and
\cite{HW} for further discussions about this problem).

Let $M$ be a II$_1$-factor, $M\subseteq B(H)$, and let $T\in M$.
If $K\subseteq H$ is a non-trivial closed $T$-invariant subspace
affiliated with $M$, and $P=P_K$ is the orthogonal projection on
$M$, then according to the decomposition, $H=K\oplus K^{\perp}$,
we can write
\begin{equation}
\label{eq9-7}
T=\begin{pmatrix} T_{11} & T_{12}\\ 0 & T_{22}\end{pmatrix},
\end{equation}
where $T_{11}=PTP$ and $T_{22}=(1-P)T(1-P)$ are elements of the
II$_1$-factors $M_1=PMP$ and $M_2=(1-P)M(1-P)$.  Let $\mu_{T_{11}}$ and
$\mu_{T_{22}}$ be the Brown measures of $T_{11}$ and $T_{22}$ computed
relative to $M_1$ and $M_2$ (respectively) then by \cite{Br}:
\begin{equation}
\label{eq9-8}
\mu_T = a\mu_{T_{11}}+(1-a)\mu_{T_{22}}
\end{equation}
where $a=\tr_M(P)$.

The main result of \cite{Haa2} is\vskip -2mm

\begin{theorem} \label{thm9-3} {\rm \cite{Haa2}}
Let $M$ be II$_1$-factor with the embedding property, and let $T\in
M$. Then for every Borel set $B\subseteq\bC$ there is a unique
$T$-invariant subspace $K$ affiliated with $M$, such that $\mu_{T_{11}}$
is concentrated on $B$ and $\mu_{T_{22}}$ is concentrated on
$\bC\backslash B$, where $T_{11}$ and $T_{22}$ are defined as in
\eqref{eq9-7}. Moreover, $\Tr_M(P_K)=\mu_T(B)$, where $P_K\in M$ is the
projection onto $K$.
\end{theorem}\vskip -2mm

\begin{remark}
\label{rem9-4}\rm  If $T$ is decomposable and $B$ is closed, then
the subspace $K$ coincide with the spectral subspace $E(B)$
characterized by \eqref{eq9-3}--\eqref{eq9-6}. However, already in
the hyperfinite II$_1$-factor $R$, there are operators $T$ which
are not decomposable.
\end{remark}\vskip -2mm

\begin{corollary} \label{cor9-5} {\rm \cite{Haa2}}
Let $T\in M$, where $M$ is a II$_1$-factor with the embedding
property. If the Brown measure $\mu_T$ of $T$ is not concentrated in a
single point, then $T$ has a non-trivial closed invariant subspace
affiliated with $M$.
\end{corollary}\vskip -2mm

\begin{remark} \label{rem9-6}\rm
Corollary \ref{cor9-5} reduced the invariant subspace problem for
II$_1$-factors $M$ with the embedding problem to operators $T\in M$ for
which $\mu_T=\delta_0$ (the Dirac-measure at 0). It can be shown that
$\mu_T=\delta_0$ if and only if
\[
\lim_{n\to\infty} ((T^*)^nT^n)^{\frac1n}=0
\]
in the strong operator topology on $M$ (cf.\ \cite{Haa2}).
\end{remark}\vskip -2mm

In the rest of this section, I will briefly outline the proof of Theorem
\ref{thm9-3}.

Let $M$ be a II$_1$-factor and let $T\in M$. Define the modified
spectrum $\sigma'(T)$ and modified spectral radius $r'(T)$ by
\begin{eqnarray*}
\sigma'(T) &=& \mbox{supp}(\mu_T)\\
r'(T) &=& \max\{ |\lambda| \mid \lambda\in\sigma'(T)\}.
\end{eqnarray*}
Then $\sigma'(T)\subseteq\sigma(T)$ and $r'(T)\le r(T)$.

The classical spectral radius formula
\begin{equation*}
r(T) = \lim_{n\to\infty} \|T^n\|^{\frac1n}
\end{equation*}
has a modified version (cf. \cite{Haa2}):
\begin{equation*}
r'(T) =
\lim_{p\to\infty}(\lim_{n\to\infty}\|T^n\|^\frac1n_{\frac{p}{n}})
\end{equation*}
where $\|S\|_p=\tr_M(|S|^p)^{\frac1p}$, $p>0$.

\vspace{.7cm}

\noindent{\bf Spectral subspace lemma 8.7}  \cite{Haa2}
\it Let $M$ be a II$_1$-factor. (Here we  do not need the embedding
property.) Let $T\in M$ and let $F\subseteq\bC$ be a closed set. Then
\begin{itemize}
\item[(a)] There exists a maximal closed $T$-invariant subspace $K$
affiliated with $M$ such that $\sigma'(T_{|K})\subseteq F$, where
$\sigma'(T_{|K})$ is the modified spectrum of the operator $T_{|K}$
considered as an element of the II$_1$-factor $P_KMP_K$ ($P_K$ is the
projection of $H$ onto $K$).
\item[(b)] Let $K(F)$ be the subspace $K$ defined by (a). Then
\[
\tr_M(P_{K(F)})\le \mu(F)
\]
for all closed subsets $F$ of $\bC$.
\end{itemize}
\rm

\vspace{.1cm}

\noindent{\bf Random distortion lemma 8.8}  \cite{Haa2}
\it Let $M$ be a II$_1$-factor with the embedding property and let $T\in
M$. Then
\begin{itemize}
\item[(a)] There exist natural numbers $k(1)<k(2)<\dots$ and $T_n\in
M_{k(n)}(\bC)$ such that
\begin{equation}
\label{eq10-5}
\sup_{n\in\bN} \|T_n\|<\infty.
\end{equation}
\item[(b)] For every non-commutative polynomial $p$ in two variables
\begin{equation}
\label{eq10-6}
\lim_{n\to\infty} \tr_{k(n)}(p(T_n,T_n^*))=\tr(p(T,T^*))
\end{equation}
where $\tr_{k(n)}$ is the normalized trace on $M_n(\bC)$.
\item[(c)] Furthermore, there exists a sequence $T'_n\in M_{k(n)}(\bC)$
such that
\begin{eqnarray}
\label{eq10-7}
\lim_{n\to\infty} \|T'_n-T_n\|_p &=&0\quad\mbox{for some $p>0$}\\
\label{eq10-8}
\lim_{n\to\infty} \Delta(T'_n-\lambda 1) &=& \Delta(T-\lambda
1)\quad\mbox{for almost all $\lambda\in\bC$}\\
\label{eq10-9}
\lim_{n\to\infty} \mu_{T'_n} &=& \mu_T\quad\mbox{weakly in Prob$(\bC)$.}
\end{eqnarray}
\end{itemize}
\rm

The embedding property is needed in (b). To pass from (b) to (c) we use
a random distortion argument where we put
\[
T'_n = T_n+\varepsilon_n X_nY^{-1}_n
\]
where $X_n,Y_n$ are random Gaussian matrices with independent entries
and $\varepsilon_n\to 0$. Subsequently Sniady proved \cite{Sn1} that by
using a different random distortion, one can obtain a stronger result,
namely in (c), \eqref{eq10-7} can be replaced by
\[
\lim_{n\to\infty}\|T'_n-T_n\|_\infty=0
\]
where $\|\cdot \|_\infty$ is the operator norm.

The random distortion lemma is used to reduce the proof of Theorem
\ref{thm9-3} to the case of $M=M_n(\bC)$ by an ultraproduct argument.
For $M=M_n(\bC)$, Theorem \ref{thm9-3} is a corollary of Jordan's normal
form.

\section{\boldmath Voiculescu's circular operator $Y$ and the strictly upper triangular
operator $T$}
\label{section 11}  \setzero\vskip-5mm \hspace{5mm }

\setcounter{equation}{0}

Prior to the proof of theorem \ref{thm9-3}, Dykema and the lecturer had
constructed invariant subspaces for special operators in factors of type
II$_1$. An example of particular interest is Voiculescu's circular
operator $Y$, which can be written as
\[
Y = \frac{1}{\sqrt{2}} (X_1 + iX_2)
\]
where $(X_1,X_2)$ is a semicircular system (cf.\ Section \ref{section
3}). The von Neumann algebra $M=VN(Y)$ generated by $Y$ is isomorphic to
$L(F_2)$ (the von Neumann associated to a free group on two
generators) which is a factor of type II$_1$. The operator $Y$ is far
from being normal and for some time it was considered a possible
counterexample for the invariant subspace problem relative to the
II$_1$-factor it generates. In \cite{HL} Larsen and the lecturer proved
that
\begin{align}
\label{eq11-1} & \mbox{$\sigma(Y)=\overline{D}$ \ (the closed unit disc
in $\bC$)}\\[7pt]
\label{eq11-2}
& \mbox{The Brown measure $\mu_Y$ of $Y$ is the uniform}  \\
\nonumber
& \mbox{distribution on $\overline{D}$, i.e. it has constant density
$\frac{1}{\pi}$.}
\end{align}\vskip -2mm

\begin{theorem}
\label{thm11-1} {\rm \cite{DH1}}
For each $r\in(0,1)$ there is a unique projection $p\in M=VN(Y)$ such
that
\begin{align}
\label{eq11-3}
& pYp=Yp\quad\mbox{(i.e. the range of $p$ is $Y$-invariant)}\\
\label{eq11-4}
& \sigma(pYp)\subseteq \{z\in\bC\mid |z|\le r\}\\
\label{eq11-5}
& \sigma((1-p)Y(1-p)) \subseteq\{z\in\bC\mid r\le |z|\le 1\}
\end{align}
where the spectra in \eqref{eq11-4} and \eqref{eq11-5} are computed
relative to $pMp$ and \newline $(1-p)M(1-p)$. Moreover
\begin{equation}
\label{eq11-6}
\tr_M(p)=r^2.
\end{equation}
\end{theorem}\vskip -2mm

This result was generalized to arbitrary $R$-diagonal elements by Sniady
and Speicher \cite{SS}. Later Dykema and the lecturer proved\vskip -2mm

\begin{theorem}
\label{thm11-2}
{\rm \cite{DH2}}
Voiculescu's circular operator is decomposable in the sense of Apostol
and Foias (see Definition \ref{def9-1}).
\end{theorem}\vskip -2mm

In \cite{DH2} we also considered the ``strictly upper triangular
operator'' $T$. It is defined in terms of its random matrix model:\vskip
-2mm

\begin{theorem/def} \label{thdef11-3} {\rm \cite{DH2}}
Let for each $n\in\bN$ \ $T_n$ denote the strictly upper triangular random
matrix
\begin{equation}
\label{eq11-7}
T_n = \begin{pmatrix} 0 & t_{11}^{(n)} & \cdots & t_{1n}^{(n)}\\
\ddots && \ddots & t_{n-1,n}^{(n)} \\
0 & && 0 \end{pmatrix}
\end{equation}
for which the entries $(t_{ij}^{(n)})_{i<j}$ are $\frac{n(n-1)}{2}$
independent identically distributed complex Gaussian random variables with
densities $\frac{n}{\pi}\exp(-n|z|^2)$, $z\in\bC$. Then there is an
operator $T$ in a II$_1$-factor $M$ such that $T_n$ converges in
*-moments to $T$, i.e.
\begin{equation}
\label{eq11-8}
\tr_M(P(T,T^*))=\lim_{n\to\infty}\bE\,\tr_n(P(T_n,T_n^*))
\end{equation}
for every non-commutative polynomial $P$. $T$ is called the strictly upper
triangular operator.
\end{theorem/def}\vskip -2mm

The strictly upper triangular operator is quasi nilpotent,
i.e. $\sigma(T)=\{0\}$, and therefore its Brown measure $\mu_T$ is equal
to $\delta_0$. In view of remark \ref{rem9-6} it could be a candidate
for a counterexample to the invariant subspace problem relative to a II$_1$-factor. However,
this is not the case:

Dykema and the lecturer proved in \cite{DH2} that
\begin{equation}
\label{eq11-8a}
\tr((T^*T)^n) = \frac{n^n}{(n+1)!},\qquad n\in\bN
\end{equation}
and in \cite{Sn2}, Sniady proved
\begin{equation}
\label{eq11-9}
\tr(((T^k)^*T^k)^n)=\frac{n^{nk}}{(nk+1)!},\qquad n,k\in\bN,
\end{equation}
a formula which was conjectured in \cite{DH2}.

Based on \eqref{eq11-9} and its proof, we recently proved\vskip -2mm

\begin{theorem}
\label{thm11-4} {\rm \cite{DH3}}
Let $T$ be as above. Put $S_k=k((T^k)^kT^k)^{\frac1k}$ and let
$F\colon [0,\pi]\to[0,1]$ be the strictly increasing function given by
$F(0)=0$, $F(\pi)=1$ and
\begin{equation}
\label{eq11-10}
F\bigg(\frac{\sin v}{v}\exp(v\cot v)\bigg) = 1-\frac{v}{\pi}+\frac{1}{\pi}
\frac{\sin^2 v}{v},\quad 0<v<\pi.
\end{equation}
Then $F(S_k)$ converges in strong operator topology to the ``diagonal
operator'' $D_0$ with matrix model
\begin{equation}
\label{eq11-11}
D_{0,n} = \begin{pmatrix} \frac1n & & & 0\\
 & \frac2n && \\
  & &  \ddots & \\
  0 &&& 1 \end{pmatrix}.
\end{equation}
In particular $D_0\in VN(T)$. Moreover $VN(T)$ is isomorphic to $L(F_2)$
and the ranges of the projections $1_{[0,t]}(D_0)$, $0<t<1$, form an
uncountable family of non-trivial invariant subspaces for $T$
affiliated with $VN(T)$.
\end{theorem}

\label{lastpage}

\end{document}